\documentclass[11pt,fleqn]{article}
\usepackage{amsfonts, epsfig, amssymb, amsmath}
\usepackage[normalem]{ulem}
\usepackage{color}
\newcommand{\ol}{\setlength{\itemsep}{0pt.}\begin{enumerate}}
\newcommand{\eol}{\end{enumerate}\setlength{\itemsep}{-\parsep}}
\newcommand{\ignore}[1]{}
\title{On $\ell_4 : \ell_2$ ratio of functions with restricted Fourier support}
\author{Naomi Kirshner and Alex Samorodnitsky}

\begin{document}
\date{}
\maketitle


\newtheorem{THEOREM}{Theorem}[section]
\newenvironment{theorem}{\begin{THEOREM} \hspace{-.85em} {\bf :}
}%
                        {\end{THEOREM}}
\newtheorem{LEMMA}[THEOREM]{Lemma}
\newenvironment{lemma}{\begin{LEMMA} \hspace{-.85em} {\bf :} }%
                      {\end{LEMMA}}
\newtheorem{COROLLARY}[THEOREM]{Corollary}
\newenvironment{corollary}{\begin{COROLLARY} \hspace{-.85em} {\bf
:} }%
                          {\end{COROLLARY}}
\newtheorem{PROPOSITION}[THEOREM]{Proposition}
\newenvironment{proposition}{\begin{PROPOSITION} \hspace{-.85em}
{\bf :} }%
                            {\end{PROPOSITION}}
\newtheorem{DEFINITION}[THEOREM]{Definition}
\newenvironment{definition}{\begin{DEFINITION} \hspace{-.85em} {\bf
:} \rm}%
                            {\end{DEFINITION}}
\newtheorem{EXAMPLE}[THEOREM]{Example}
\newenvironment{example}{\begin{EXAMPLE} \hspace{-.85em} {\bf :}
\rm}%
                            {\end{EXAMPLE}}
\newtheorem{CONJECTURE}[THEOREM]{Conjecture}
\newenvironment{conjecture}{\begin{CONJECTURE} \hspace{-.85em}
{\bf :} \rm}%
                            {\end{CONJECTURE}}
\newtheorem{MAINCONJECTURE}[THEOREM]{Main Conjecture}
\newenvironment{mainconjecture}{\begin{MAINCONJECTURE} \hspace{-.85em}
{\bf :} \rm}%
                            {\end{MAINCONJECTURE}}
\newtheorem{PROBLEM}[THEOREM]{Problem}
\newenvironment{problem}{\begin{PROBLEM} \hspace{-.85em} {\bf :}
\rm}%
                            {\end{PROBLEM}}
\newtheorem{QUESTION}[THEOREM]{Question}
\newenvironment{question}{\begin{QUESTION} \hspace{-.85em} {\bf :}
\rm}%
                            {\end{QUESTION}}
\newtheorem{REMARK}[THEOREM]{Remark}
\newenvironment{remark}{\begin{REMARK} \hspace{-.85em} {\bf :}
\rm}%
                            {\end{REMARK}}

\newcommand{\thm}{\begin{theorem}}
\newcommand{\lem}{\begin{lemma}}
\newcommand{\pro}{\begin{proposition}}
\newcommand{\dfn}{\begin{definition}}
\newcommand{\rem}{\begin{remark}}
\newcommand{\xam}{\begin{example}}
\newcommand{\cnj}{\begin{conjecture}}
\newcommand{\mcnj}{\begin{mainconjecture}}
\newcommand{\prb}{\begin{problem}}
\newcommand{\que}{\begin{question}}
\newcommand{\cor}{\begin{corollary}}
\newcommand{\prf}{\noindent{\bf Proof:} }
\newcommand{\ethm}{\end{theorem}}
\newcommand{\elem}{\end{lemma}}
\newcommand{\epro}{\end{proposition}}
\newcommand{\edfn}{\bbox\end{definition}}
\newcommand{\erem}{\bbox\end{remark}}
\newcommand{\exam}{\bbox\end{example}}
\newcommand{\ecnj}{\bbox\end{conjecture}}
\newcommand{\emcnj}{\bbox\end{mainconjecture}}
\newcommand{\eprb}{\bbox\end{problem}}
\newcommand{\eque}{\bbox\end{question}}
\newcommand{\ecor}{\end{corollary}}
\newcommand{\eprf}{\bbox}
\newcommand{\beqn}{\begin{equation}}
\newcommand{\eeqn}{\end{equation}}
\newcommand{\wbox}{\mbox{$\sqcap$\llap{$\sqcup$}}}
\newcommand{\bbox}{\vrule height7pt width4pt depth1pt}
\newcommand{\qed}{\bbox}
\def\sup{^}

\def\H{\{0,1\}^n}

\def\S{S(n,w)}

\def\g{g_{\ast}}
\def\xop{x_{\ast}}
\def\y{y^{\ast}}
\def\z{z_{\ast}}

\def\f{\tilde f}

\def\n{\lfloor \frac n2 \rfloor}

\def \E{\mathop{{}\mathbb E}}
\def \R{\mathbb R}
\def \Z{\mathbb Z}
\def \F{\mathbb F}
\def \S{\mathbb S}

\def \x{\textcolor{red}{x}}
\def \r{\textcolor{red}{r}}
\def \Rc{\textcolor{red}{R}}

\def \noi{{\noindent}}

\def \iff{~~~~\Leftrightarrow~~~~}

\def \queq {\quad = \quad}

\def\<{\left<}
\def\>{\right>}
\def \({\left(}
\def \){\right)}

\def \e{\epsilon}
\def \l{\lambda}

\def\Tp{Tchebyshef polynomial}
\def\Tps{TchebysDeto be the maximafine $A(n,d)$ l size of a code with distance $d$hef polynomials}
\newcommand{\rarrow}{\rightarrow}

\newcommand{\larrow}{\leftarrow}

\overfullrule=0pt
\def\setof#1{\lbrace #1 \rbrace}

\begin{abstract}
Given a subset $A \subseteq \H$, let $\mu(A)$ be the maximal ratio between $\ell_4$ and $\ell_2$ norms of a function whose Fourier support is a subset of $A$.\footnote{Strictly speaking, we consider the fourth power of this ratio, since it is easier to work with.} We make some simple observations about the connections between $\mu(A)$ and the additive properties of $A$ on one hand, and between $\mu(A)$ and the uncertainty principle for $A$ on the other hand. One application obtained by combining these observations with results in additive number theory is a stability result for the uncertainty principle on the discrete cube.

Our more technical contribution is determining $\mu(A)$ rather precisely, when $A$ is a Hamming sphere $S(n,k)$ for all $0 \le k \le n$.

\end{abstract}

\section{Introduction}

\noi Let $A$ be a subset of the discrete cube $\H$. Consider the subspace $V = V(A)$ of functions on $\H$ whose Fourier support is a subset of $A$. That is, for any function $f \in V$, the expansion of $f$ in terms of the Walsh-Fourier characters $f = \sum_{\alpha} \widehat{f}(\alpha) W_{\alpha}$ is supported on $\alpha \in A$. Let
\[
\mu(A) \queq \max_{f \in V, f \not = 0} \(\frac{\|f\|_4}{\|f\|_2}\)^4 \queq \max_{f \in V, f \not = 0} \frac{\E f^4}{\E^2 f^2},
\]
where the expectation on the RHS is w.r.t. the uniform measure on $\H$.

\noi The quantity $\mu(A)$ is well-investigated, especially when $A$ is a Hamming ball or a Hamming sphere, since in this case it is closely related to the hypercontractive property of the noise operator on the discrete cube. In particular, it is know that for a Hamming ball of radius $k$, we have $\mu(A) \le 9^k$ \cite{Bonami}, and for a Hamming sphere of radius $k$, for a slowly growing $k$, we have $\mu(A) = \Theta\(9^k / \sqrt{k}\)$ \cite{O'Donnel}.

\noi We make several simple observations, connecting between $\mu(A)$ and the additive properties of $A$ on one hand, and between $\mu(A)$ and the uncertainty principle for functions in $V(A)$ on the other hand. Connections of this kind have already been explored in \cite{KM, PS} (between $\mu(A)$, or closely related quantities, to the uncertainty principle), and by \cite{Hastad-personal} (between $\mu(A)$ and the additive properties of $A$).

\noi {\bf Additive structure of $A$}: For $x \in A + A$, let $M_x = \{(a,b) \in A \times A: a + b = x\}$.  Let $m(A) = 1 + \max_{x \not = 0} |M_x|$. Thus $m(A)$ is the maximal multiplicity of a non-zero element in $A+A$ (plus one). Let  $E_2(A,A)$ be the {\it additive energy} of $A$ \cite{TV}. This is the number of $4$-cycles in $A$:
\[
E_2(A,A) \queq \Big | \{(a,b,c,d) \in A^4: a+b+c+d = 0\} \Big |
\]

\noi We observe that

\pro
\label{pro:additive}
For any subset $A \subseteq \H$ holds
\begin{enumerate}
\item
$
\mu(A) ~ \le ~ |A|
$
\item
$
\mu(A) ~\le~ m(A)
$
\item
$
\max_{B \subseteq A} \frac{E_2(B,B)}{|B|^2} ~\le~ \mu(A) ~ \le ~ O\(\log^3(|A|)\) \cdot \max_{B \subseteq A} \frac{E_2(B,B)}{|B|^2}
$
\end{enumerate}
\epro

\noi Let us describe one application of this proposition. The following result has been proved in  \cite{GT}. Let $A$ be a Hamming ball of radius $k$, and let $B$ and $C$ be subsets of $A$. Then $|B+C| \ge \frac{|B||C|}{9^k}$. We rederive this result as follows. Recall that $\mu(A) \le 9^k$. Hence we have
\beqn
\label{gt-ineq}
|B+C| ~\ge~ \frac{|B|^2 |C|^2}{\sqrt{E_2(B,B)}  \sqrt{E_2(C,C)}} ~=~ \frac{|B| |C|}{\sqrt{\frac{E_2(B,B)}{|B|^2}} \sqrt{\frac{E_2(C,C)}{|C|^2}}} ~\ge~\frac{|B| |C|}{\mu(A)} ~\ge~ \frac{|B||C|}{9^k}.
\eeqn

\noi For the first inequality (which follows by a simple application of the Cauchy-Schwarz inequality) see e.g., \cite{TV}. The second inequality follows from the third claim of the proposition.

\rem
\label{rem:eigenvalue}
Another way to obtain $\mu(A)$ is as the maximal eigenvalue of a certain symmetric $|A| \times |A|$ matrix. Such matrices and their relevance to the additive structure of $A$ were considered in \cite{Shkredov}. Specifically, denoting by $\l(M)$ the maximal eigenvalue of a matrix $M$, it is not hard to see that
\[
\mu(A) \queq \max_{y:A \rarrow \R, \|y\|_2 = 1}  \l\(T^{y \circ y}_A\),
\]
where (in the notation of \cite{Shkredov}) $T = T^{y \circ y}_A$ is the $A \times A$ matrix with rows and columns indexed by the elements of $A$, such that $T\(a_1,a_2\) = \sum_{\(b_1,b_2\) \in M_{a_1 + a_2}} y\(a_1\) \cdot y\(a_2\)$.
\erem

\noi {\bf Uncertainty principle}: The uncertainty principle for the discrete cube (see e.g., \cite{Don}) states that for a non-zero function $f$ on $\H$ holds:
\beqn
\label{ineq-uncertain}
|supp(f)| \quad \ge \quad \frac{2^n}{|supp\(\widehat{f}\)|}
\eeqn

\noi The following claim is an immediate consequence of the Cauchy-Schwarz inequality.
\lem
\label{lem:uncertain}
For a non-zero function $f$ on $\H$, let let $A = supp\(\widehat{f}\)$. Then
\[
|supp(f)| \quad \ge \quad \frac{2^n}{\mu(A)}
\]
\elem

\noi This strengthens (\ref{ineq-uncertain}), by the first claim of Proposition~\ref{pro:additive}.

\noi A quantitative version of (\ref{ineq-uncertain}) was proved in \cite{PS}. For any $0 < \delta < 1$, there exists an $\e > 0$ depending on $\delta$, such that for any two subsets $A$, $B$ of $\H$ with $|A| \cdot |B| \le 2^{(1-\delta) n}$ holds: if $f$ is a non-zero function with $\widehat{f}$ supported on $A$, then
$\frac{1}{2^n} \sum_{b \in B} f^2(b) \le (1 - \e) \cdot \|f\|^2_2$.

\noi The following claim is a strengthening of this result.
\lem
\label{lem:str-PS}
Let $0 < \delta < 1$. Then for any two subsets $A$, $B$ of $\H$ with $\mu(A) \cdot |B| \le 2^{(1-\delta) n}$ holds: if $f$ is a non-zero function with $\widehat{f}$ supported on $A$, then
\[
\frac{1}{2^n} \sum_{b \in B} f^2(b) \quad \le \quad 2^{- \frac{\delta n}{2}} \cdot \|f\|^2_2
\]
\elem

\noi Combining Lemma~\ref{lem:uncertain} with Proposition~\ref{pro:additive} gives the following corollary.
\cor
\label{cor:uncertain}
For a non-zero function $f$ on $\H$, let $A = supp\(\widehat{f}\)$. Then
\[
|supp(f)| \quad \ge \quad \frac{2^n}{m(A)} \quad \text{and} \quad |supp(f)| \quad \ge \quad \Omega\(\frac{1}{\log^3(|A|)}\) \cdot \frac{2^n}{\max_{B \subseteq A} \frac{E_2(B,B)}{|B|^2}}
\]
\ecor

\noi Up to negligible factors, both inequalities strengthen (\ref{ineq-uncertain}), since $m(A) \le |A| + 1$, and $E_2(B)~\le~|B|^3$.

\noi Combining the second inequality in Corollary~\ref{cor:uncertain} with results from additive number theory describing the structure of sets with large energy and small doubling (\cite{BGS}, \cite{Sanders}) leads to a stability version of (\ref{ineq-uncertain}). It is known that (\ref{ineq-uncertain}) holds with equality if and only if $\widehat{f}$ is a characteristic function of an affine subspace of $\H$. We show that even if equality is replaced with 'near equality', the support of $\widehat{f}$ will be similar to a linear subspace, in the appropriate sense. {\it Notation}: let $\<B\>$ denote the linear span of a subset $B \subseteq \H$.

\pro
\label{pro:stab-uncert}
Let $f$ be a non-zero function on $\H$ with $|supp(f)| \cdot |supp\(\widehat{f}\)| \le C \cdot 2^n$. Let $A = supp\(\widehat{f}\)$. Let $C' = C \cdot \log(|A|)$. There exists a subset $A' \subseteq A$ such that:
\begin{itemize}
\item
\[
|A'| \quad \ge \quad C'^{-O\(\log^3 C'\)} \cdot |A|
\]
and
\item
\[
|\<A'\>| \quad \le \quad  |A|,
\]
\end{itemize}
with asymptotic notation hiding absolute constants.
\epro

\noi The third claim of Proposition~\ref{pro:additive} leads to the following natural question: which sets $A \subseteq \H$ have the 'hereditary' property $\frac{E_2(A,A)}{|A|^2} \ge \frac{E_2(B,B)}{|B|^2}$, for all subsets $B \subseteq A$. It is easy to see that this holds if $A$ is a subspace. We show that, up to lower order terms, this is also true for a Hamming sphere. We distinguish between two cases: the radius of the sphere is small compared to $n$, or the radius of the sphere is allowed to grow arbitrarily in $n$. For the first case, we have the following proposition:
\pro
\label{pro:small sphere}
Let $A = S(n,k)$ be a Hamming sphere of radius $k$, for $k = o(\sqrt{n})$. Then
\[
\mu(A) \quad \le \quad \Big(1 + o_n(1)\Big) \cdot \frac{E_2(A,A)}{|A|^2}
\]
\epro

\noi For general $k$ we have the following result, which is the most technical part of this paper.

\thm
\label{thm:sphere}
Let $A = S(n,k)$ be a Hamming sphere of radius $k$, for $0 \le k \le n/2$. Let $r(x) = \frac{3 - \sqrt{1 + 8(1-2x)^2}}{8}$, and let $\psi$ be a function on $\left[0,\frac 12\right]$ defined by
\[
\psi(x) \queq H\big(2r(x)\big) + 4r(x) + 2\big(1-2r(x)\big) \cdot H\(\frac{x-r(x)}{1-2r(x)}\) - 2H(x).
\]

\noi Then
\begin{enumerate}
\item
\[
\mu(A) \quad \le \quad 2^{n \psi\(\frac{k}{n}\)}
\]
\item
\[
2^{n \psi\(\frac{k}{n}\)} \quad \le \quad O\(k^{3/2}\) \cdot \frac{E_2(A,A)}{|A|^2}
\]
\end{enumerate}
\ethm

\noi In the light of these results it is natural to make the following conjecture.
\cnj
\label{cnj:sphere}
Let $A = S(n,k)$ be a Hamming sphere of radius $k$, for $1 \le k \le n/2$. Then
\[
\mu(A) \queq \frac{E_2(A,A)}{|A|^2}
\]
\ecnj

\rem
\label{rem:Krawchouk}
Proposition~\ref{pro:small sphere} and Theorem~\ref{thm:sphere} show that among all homogeneous polynomials $f$ of degree $k$ the maximum of the ratio $\frac{\|f\|_4}{\|f\|_2}$ is (essentially) attained for the sum of all weight $k$ monomials (the $k^{th}$ {\it Krawchouk polynomial} $K_k$). Equivalently, these results essentially determine the $\|\cdot\|_{2 \rarrow 4}$ norm of the projection operator $P_k:~f \rarrow \sum_{|\alpha| = k} \widehat{f}(\alpha) W_{\alpha}$ (see \cite{P2} where the norms of these operators are investigated).
We refer to \cite{NO} (and the references therein) and to Sections 4.2 and 4.3 in \cite{PS} for other results in this direction.

\erem

\noi Theorem~\ref{thm:sphere} can be applied to extend the result of \cite{GT} (whose alternative derivation was given in (\ref{gt-ineq})) to larger values of $k$. We start with observing that a simple modification of the proof of the theorem shows its bound to hold for Hamming balls as well.

\cor
\label{cor:ball}
Let $A = B(n,k)$ be a Hamming ball of radius $k$, for $0 \le k \le n/2$. Then:
\begin{enumerate}

\item
$\mu(A) ~\le~ 2^{n \psi\(\frac{k}{n}\)}$.

\item
$2^{n \psi\(\frac{k}{n}\)} ~\le~ \min \left\{9^k, 2^n\right\}$,

\noi with equality only at $k = 0$, where the LHS is $1$ and at $k = n/2$, where the LHS is $2^n$.

\end{enumerate}

\ecor

\noi The second claim of the corollary shows it to extend the bound $\mu(A) \le 9^k$ (\cite{Bonami}). Using its first claim in (\ref{gt-ineq}) leads to the following result (which we state slightly more generally).

\cor
\label{gt-stronger}
Let $B$ be a subset of a Hamming ball of radius $k_1$, and let $C$ be a subset of a Hamming ball of radius $k_2$. Then
\[
|B+C| \quad \ge \quad \frac{|B||C|}{2^{\frac{n}{2} \cdot \(\psi\(\frac{k_1}{n}\)+\psi\(\frac{k_2}{n}\)\)}}.
\]

\ecor

\noi This paper is organized as follows. We prove Proposition~\ref{pro:small sphere} in Section~\ref{sec:small sphere}, and Theorem~\ref{thm:sphere} with Corollary~\ref{cor:ball} in Section~\ref{sec:sphere}. All the remaining claims are proved in Section~\ref{sec:simple}.

\section{Simple proofs}
\label{sec:simple}

\noi In this section we prove all the observations stated in the introduction, except for Proposition~\ref{pro:small sphere} and Theorem~\ref{thm:sphere}.

\noi Our starting point is the following characterization of $\mu(A)$. Let $\S = {\S}^{|A|-1}$ denote the Euclidean sphere of dimension $|A|-1$. We will assume the vectors in $\S$ to be indexed by elements of $A$ (in other words, a vector $y \in \S$ is a function from $A$ to $\R$, with unit $\ell_2$ norm). Then $\mu(A)$ is the maximal value of the following real valued function on $\S$ (recall that $M_x = \{(a,b) \in A \times A: a + b = x\}$):
\beqn
\label{mu-char}
\mu(A) \queq \max_{y \in \S} F(y), \quad \mbox{where} \quad F(y) \queq \sum_{x \in A+A} \(\sum_{(a,b) \in M_x} y_a y_b\)^2.
\eeqn

\noi To see this, note that each $y \in \S$ represents a Fourier expansion of a function $f = \sum_{a \in A} y_a W_a$ of $\ell_2$ norm $1$, and $F(y) = \E f^4 = \|f\|^4_4$.

\subsubsection*{Proof of Proposition~\ref{pro:additive}}

\noi We start with the first claim of the proposition. Applying the Cauchy-Schwarz inequality, for any $y \in S^{|A|-1}$ holds
\[
F(y) \queq \sum_{x \in A+A} \(\sum_{(a,b) \in M_x} y_a y_b\)^2 \quad \le \quad \sum_{x \in A+A} |M_x| \cdot \sum_{(a,b) \in M_x} y^2_a y^2_b \quad \le
\]
\[
\(\max_{x \in A+A} |M_x|\) \cdot \sum_{x \in A+A} \sum_{(a,b) \in M_x} y^2_a y^2_b \queq  \max_{x \in A+A} |M_x| \queq |M_0| \queq |A|,
\]
completing the proof.

\noi The second claim is proved similarly. For any $y \in \S$ holds
\[
F(y) \queq \sum_{x \in A+A} \(\sum_{(a,b) \in M_x} y_a y_b\)^2 \quad \le \quad 1 + \sum_{x \in A+A \setminus \{0\}} |M_x| \cdot \sum_{(a,b) \in M_x} y^2_a y^2_b \quad \le
\]
\[
1 + \(\max_{x \in A+A \setminus \{0\}} |M_x|\) \cdot \sum_{x \in A+A} \sum_{(a,b) \in M_x} y^2_a y^2_b \queq  1 + \max_{x \in A+A \setminus \{0\}} |M_x|,
\]

\noi We continue to the third claim, starting with the lower bound. Note that for any subset $B \subseteq A$ holds $F\(\frac{1_B}{\sqrt{|B|}}\) = \frac{E_2(B,B)}{|B|^2}$. Hence, by (\ref{mu-char}),
\[
\mu(A) \quad \ge \quad \max_{B \subseteq A} F\(\frac{1_B}{\sqrt{|B|}}\) \queq \max_{B \subseteq A} \frac{E_2(B,B)}{|B|^2}.
\]

\noi We pass to the upper bound on $\mu(A)$. Let $\y \in \S$ such that $F(\y) = \mu(A)$. We may assume, w.l.o.g, that the vector $\y$ is nonnegative. Let $f = \sum_{a \in A} \y_a W_a$. Then $\E f^4 = \mu(A)$.

\noi We introduce some notation: For $i \ge 1$, let $A_i = \left\{a \in A: ~2^{-i} < \y_a \le 2^{-(i-1)}\right\}$. Let $f_i = \sum_{a \in A_i} \y_a W_a$. Let $h_i = \sum_{a \in A_i} W_a$. Finally, let $N = \lceil\frac12 \log_2(|A|)\rceil + 2$.

\noi We have $f = \sum_{i=1}^{\infty} f_i$, where the summation on the RHS is, of course, finite. Let $k = \sum_{i=N+1}^{\infty} f_i$. Then $k = \sum_{a \in A} z_a W_a$, with $|z_a| \le 2^{-(N-1)} \le \frac{1}{2\sqrt{|A|}}$ for all $a \in A$. Hence $\sum_{a \in A} z^2_a \le \frac{1}{4}$ and therefore, by the $4$-homogeneity of $F$, we get
$\E k^4 = F(z) \le \frac{\mu(A)}{16}$.

\noi Let $t = f-k = \sum_{i=1}^N f_i$. By the convexity of the function $x^4$ and by Jensen's inequality, we have $\E f^4 = \E (k+t)^4 \le 8 \cdot \(\E k^4 + \E t^4\)$. It follows that $\E t^4 \ge \frac18 \cdot \E f^4 - \frac{\mu(A)}{16} \ge \frac{1}{16} \cdot \E f^4$. So, to prove the claim it suffices to upperbound $\E t^4$, which we proceed to do.

\noi By Jensen's inequality,
$\E t^4 = \E \(\sum_{i=1}^N f_i\)^4 \le N^3 \cdot \sum_{i=1}^N \E f^4_i$. For $1 \le i \le N$, let $y^{(i)} = 1_{A_i} \cdot \y$. Then $\E f^4_i = F\(y^{(i)}\) \le 2^{-4(i-1)} \cdot F\(1_{A_i}\) = 2^{-4(i-1)} \cdot\E h^4_i$.
Hence,
\[
\E t^4 \quad \le \quad N^3 \cdot \sum_{i=1}^N 2^{-4(i-1)} \E h^4_i \queq 16 N^3 \cdot \sum_{i=1}^N 2^{-4i} \E h^4_i
\]
Next, observe that the functions $\{f_i\}_{i=1}^N$ are orthogonal, and hence
\[
\sum_{i=1}^N 2^{-2i} |A_i| \queq \sum_{i=1}^N 2^{-2i} \E  h^2_i \quad \le \quad \sum_{i=1}^N \E f^2_i \queq \E t^2 \quad \le \quad 1.
\]

\noi It follows that
\[
\sum_{i=1}^N 2^{-4i} \E h^4_i \queq \sum_{i=1}^N 2^{-4i} E_2\(A_i, A_i\) \queq \sum_{i=1}^N \(2^{-4i} |A_i|^2\) \cdot \frac{E_2\(A_i, A_i\)}{|A_i|^2} \quad \le \quad
\]
\[
\max_{1 \le i \le N} \frac{E_2\(A_i, A_i\)}{|A_i|^2} \cdot \sum_{i=1}^N 2^{-4i} |A_i|^2 \quad \le \quad \max_{1 \le i \le N} \frac{E_2\(A_i, A_i\)}{|A_i|^2} \cdot \sum_{i=1}^N 2^{-2i} |A_i| \quad \le \quad \max_{1 \le i \le N} \frac{E_2\(A_i, A_i\)}{|A_i|^2}
\]
And hence, recalling that $N = O(\log |A|)$,
\[
\E t^4 \quad \le \quad  16 N^3 \cdot \sum_{i=1}^N 2^{-4i} \E h^4_i \quad \le \quad O\(\log^3(|A|)\) \cdot\max_{B \subseteq A} \frac{E_2(B,B)}{|B|^2},
\]
concluding the proof of the upper bound and of the proposition.

\eprf

\subsubsection*{Proof of Lemma~\ref{lem:uncertain}}

\noi Let $f$ be a non-zero function on $\H$, with $A = supp\(\widehat{f}\)$. Let $B = supp(f)$. Then, by the Cauchy-Schwarz inequality,
\[
{\E}^2 f^2 \queq {\E}^2 f^2 \cdot 1_B \quad \le \quad \E f^4 \cdot \E 1_B.
\]
Hence, by the definition of $\mu(A)$,
\[
|B| \queq 2^n \cdot \E 1_B \quad \ge \quad 2^n \cdot \frac{\E^2 f^2}{\E f^4} \quad \ge \quad \frac{2^n}{\mu(A)}.
\]

\eprf

\subsubsection*{Proof of Lemma~\ref{lem:str-PS}}

\noi Let $f$ be a non-zero function on $\H$, with $A = supp\(\widehat{f}\)$. Let $B \subseteq \H$ satisfy $|A| \cdot |B| = 2^{(1-\delta) n}$. Let $\frac{1}{2^n} \sum_{b \in B} f^2(b) = c \cdot \|f\|^2_2 = c \cdot \E f^2$. Then, by the Cauchy-Schwarz inequality,
\[
c^2 \cdot {\E}^2 f^2 \queq \(\frac{1}{2^n} \sum_{b \in B} f^2(b)\)^2 \queq {\E}^2 f^2 \cdot 1_B \quad \le \quad \E f^4 \cdot \E 1_B.
\]
Hence, by the definition of $\mu(A)$,
\[
c^2 \quad \le \quad \frac{\E f^4}{\E^2 f^2} \cdot \frac{|B|}{2^n} \quad \le \quad \frac{\mu(A) \cdot |B|}{2^n} \queq 2^{-\delta n}.
\]

\eprf

\subsubsection*{Proof of Proposition~\ref{pro:stab-uncert}}

\noi Let $f$ be a non-zero function on $\H$ with $|supp(f)| \cdot |supp\(\widehat{f}\)| \le C \cdot 2^n$. Let $A  = supp\(\widehat{f}\)$, and let $B \subseteq A$ be the subset of $A$ for which the ratio $\frac{E_2(B,B)}{|B|^2}$ is maximal. By the third claim of Proposition~\ref{pro:additive}, we have
\[
\frac{C \cdot 2^n}{|A|} \quad \ge \quad |supp(f)| \quad \ge \quad \Omega\(\frac{1}{\log^3(|A|)}\) \cdot \frac{2^n}{\frac{E_2(B,B)}{|B|^2}}.
\]
Rearranging, this gives $\frac{E_2(B,B)}{|B|^2} \ge \frac{1}{C} \cdot \Omega\(\frac{|A|}{\log^3(|A|)}\)$. Since $E_2(B,B) \le |B|^3$, this implies $|B| \ge \frac{1}{C} \cdot \Omega\(\frac{|A|}{\log^3(|A|)}\)$. Hence
\beqn
\label{A-large-energy}
E_2(A,A) \quad \ge \quad E_2(B,B) \quad \ge \quad  \frac{1}{C} \cdot \Omega\(\frac{|A||B|^2}{\log^3(|A|)}\) \quad \ge \quad  \frac{1}{C^3} \cdot \Omega\(\frac{|A|^3}{\log^9(|A|)}\)
\eeqn

\noi We quote two results from additive number theory (without stating the best known values of various constants):
\begin{itemize}

\item
\cite{BGS}: Let $A \subseteq \H$ with $E_2(A,A) \ge c \cdot |A|^3$. Then there is a subset $A_1 \subseteq A$ with $|A_1| \ge \Omega\(c^{\Theta(1)}\) \cdot |A|$ and $|A_1 + A_1| \le O\(c^{-\Theta(1)}\) \cdot |A_1|$.

\item
\cite{Sanders}: Let $A_1 \subseteq \H$ with $|A_1 + A_1| \le c_1 \cdot |A_1|$. Then there is a subset $A' \subseteq A_1$ with $|A'| \ge c_1^{-O\(\log^3 c_1\)} \cdot |A_1|$ and $|\<A'\>| \le |A_1|$.

\end{itemize}

\noi The claim of the proposition follows by combining these two results with (\ref{A-large-energy}).

\eprf

\section{Proof of Proposition~\ref{pro:small sphere}}
\label{sec:small sphere}

\noi Let $A = S(n,k)$. Then $A + A = S(n,0) \cup S(n,2) \cup ... \cup S(n,2k)$. We partition the function $F$ in (\ref{mu-char}) as $F = \sum_{t = 0}^k F_t$, where
\[
F_t \queq \sum_{x \in S(n,2t)} \(\sum_{(a,b) \in M_x} y_a y_b\)^2.
\]

\noi Clearly $F_0 \equiv 1$. We claim that for any $1 \le t \le k$ and for any $y \in \S$ holds
\beqn
\label{S_t}
F_t(y) \quad \le \quad {{2t} \choose t} \cdot {k \choose t}^2
\eeqn

\noi To see this, let $1 \le t \le k$ and let $x \in S(n,2t)$. Consider a representation $x = u + v$ with $u, v \in S(n,k)$. Note that each such representation corresponds to a partition of $x$ into two parts $x_1$ and $x_2$ of weight $t$ each, and a choice of an additional vector $w$ of weight $k-t$ disjoint from $x$, such that, slightly informally, $u = x_1 w$ and $v = x_2 w$ (that is $u$ is a concatenation of $x_1$ and $w$ and similarly for $v$).

\noi Let us denote the set of the ${{2t} \choose t}$ partitions of $x$ into two halves $x_1$ and $x_2$ by $P(x)$. Each partition $\alpha = \(x_1,x_2\) \in P(x)$ defines a subsum $s_{\alpha} = \sum_{w} y_{x_1 w} y_{x_2 w}$ of $s_x := \sum_{(a,b) \in M_x} y_a y_b$.
Clearly $s_x = \sum_{\alpha \in P(x)} s_{\alpha}$. By the Cauchy-Schwarz inequality, $s^2_x \le |P(x)|~\cdot~\sum_{\alpha \in P(x)} s^2_{\alpha} = {{2t} \choose t}~\cdot~\sum_{\alpha \in P(x)} s^2_{\alpha}$. Summing up, we have
\[
F_t(y) \queq \sum_{x \in S(n,2t)} s^2_x \quad \le \quad {{2t} \choose t} \cdot \sum_{x \in S(n,2t)} \sum_{\alpha \in P(x)} s^2_{\alpha}.
\]
Hence, (\ref{S_t}) will be implied by the following lemma.

\lem
\label{lem:St}
For any $y \in \S$ holds
\[
\sum_{x \in S(n,2t)} \sum_{\alpha \in P(x)} s^2_{\alpha}(y) \quad \le \quad {k \choose t}^2.
\]
\elem

\prf (Of the lemma)

\noi We apply the Cauchy-Schwarz inequality to bound each of the summands. For $x \in S(n,2t)$ and $\alpha = \(x_1,x_2\) \in P(x)$ we have
\[
s^2_{\alpha}(y) \queq \(\sum_{w} y_{x_1 w} y_{x_2 w}\)^2 \quad \le \quad \(\sum_w y^2_{x_1 w}\) \cdot \(\sum_w y^2_{x_2 w}\) \queq \sum_{w_1,w_2} y^2_{x_1 w_1} \cdot y^2_{x_2 w_2}
\]
That is,
\[
\sum_{x \in S(n,2t)} \sum_{\alpha \in P(x)} S^2_{\alpha}(y) \quad \le \quad \sum_{x \in S(n,2t)} ~\sum_{\(x_1,x_2\) \in P(x)} ~\sum_{w_1,w_2} y^2_{x_1 w_1} \cdot y^2_{x_2 w_2},
\]
where the inner sum goes over all $(k-t)$-bit strings $w_1, w_2$ disjoint with $x$.

\noi We will argue that for any two elements $a$ and $b$ of $S(n,k)$, the product $y^2_a y^2_b$ appears on the RHS at most ${k \choose t}^2$ times, and hence the RHS is bounded from above by ${k \choose t}^2  \cdot \sum_{a,b \in S(n,k)} y^2_a y^2_b = {k \choose t}^2  \cdot \(\sum_{a \in S(n,k)} y^2_a\)^2 = {k \choose t}^2$.

\noi In fact, given $a$ and $b$, there are at most ${k \choose t}$ ways to choose a $t$-subset $x_1 \subseteq a$ of $a$, and at most ${k \choose t}$ ways to choose a $t$-subset $x_2 \subseteq b$ of $b$. After choosing $\{x_i\}$, their complements $\{w_j\}$ are determined uniquely by $\{x_i\}$, $a$ and $b$.
\eprf

\noi This completes the proof of (\ref{S_t}). Summing up over $t$, we get
\[
\mu(A) \quad \le \quad \sum_{t=0}^k \max_{y \in \S} F_t(y) \quad \le \quad \sum_{t=0}^k {{2t} \choose t} \cdot {k \choose t}^2.
\]

\noi We proceed to compare this bound to $\frac{E_2(A,A)}{|A|^2}$. We have
\[
\frac{E_2(A,A)}{|A|^2} =  \frac{1}{|A|^2} \cdot \sum_{x \in A+A} |M_x|^2 =
\frac{1}{{n \choose k}^2} \cdot \sum_{t=0}^k \sum_{x \in S(n,2t)} | M_x |^2 = \frac{1}{{n \choose k}^2} \cdot \sum_{t = 0}^k {n \choose {2t}} \({{2t} \choose t} \cdot {{n-2t} \choose {k-t}}\)^2
\]

\noi It is easy to see that for $r = o(n)$ holds $\(1-o_n(1)\) \cdot e^{-r^2/n} \cdot n^r\le \frac{n!}{(n-r)!} \le n^r$. This implies (following a simple calculation) that for $k = o(n)$ we can lowerbound $\frac{E_2(A,A)}{|A|^2}$ by :
\beqn
\label{small-k}
\(1-o_n(1)\) e^{\frac{-2k^2}{n}} \cdot \frac{\(k!\)^2}{n^{2k}} \cdot \sum_{t=0}^k \frac{n^{2t}}{(2t)!} \cdot {{2t} \choose t}^2 \cdot \frac{n^{2k-2t}}{\((k-t)!\)^2} \queq
\(1-o_n(1)\) e^{\frac{-2k^2}{n}} \cdot \sum_{t=0}^k {{2t} \choose t} \cdot {k \choose t}^2
\eeqn

\noi Taking $k = o(\sqrt{n})$, this implies $\mu(A) \le (1 + o_n(1)) \cdot \frac{E_2(A,A)}{|A|^2}$, completing the proof of the proposition.

\section{Proof of Theorem~\ref{thm:sphere}}
\label{sec:sphere}

\noi Let $A = S(n,k)$. It will be convenient to use a notation which makes explicit the dependence of $\mu(A)$ and $\frac{E_2(A,A)}{|A|^2}$ on the parameters $n$ and $k$. We let $R(n,k) = \mu(A)$ and $r(n,k) =  \frac{E_2(A,A)}{|A|^2}$. Recall (see Section~\ref{sec:small sphere}) that
\[
r(n,k) \queq \frac{1}{{n \choose k}^2} \cdot \sum_{t = 0}^k {n \choose {2t}} \({{2t} \choose t} \cdot {{n-2t} \choose {k-t}}\)^2
\]

\noi We also let $s_t(n,k) = \frac{1}{{n \choose k}^2} \cdot {n \choose {2t}} \({{2t} \choose t} \cdot {{n-2t} \choose {k-t}}\)^2$. Thus $r(n,k) = \sum_{t=0}^k s_t(n,k)$.

\noi We start with the first (and main) claim of Theorem~\ref{thm:sphere} and rewrite it in this notation.
\beqn
\label{thm-claim}
R(n,k) \quad \le \quad 2^{n \psi\(\frac{k}{n}\)}.
\eeqn

\noi The main step in the proof of (\ref{thm-claim}) is the following somewhat weaker claim.
\pro
\label{pro:shere-main}
There exists an absolute constant $C > 0$ so that for all $1 \le k \le n/2$ holds
\[
R(n,k) \quad \le \quad C \cdot 2^{5 n/\log(n)} \cdot r(n,k).
\]
\epro

\noi We will also need the following technical lemma. From now on, all logarithms are to base $2$. Let $t_1(n,k) = \frac{3n - \sqrt{n^2 + 8(n-k)^2}}{8}$.
\lem
\label{lem:max-r}
Let $n$ be sufficiently large, and let $\frac{n}{\log n} \le k \le \frac n2 - \frac{n}{\log n}$. Then
\[
\max_{0 \le t \le k} s_t(n,k) \queq \max_{t \in t_1(n,k) \pm \sqrt{n \log n}} ~s_t(n,k).
\]
\elem

\noi We will prove Proposition~\ref{pro:shere-main} and Lemma~\ref{lem:max-r} below. First we show how they imply (\ref{thm-claim}). It will be convenient to work with the following modification of the function $\psi$. It Let $\phi$ be a function on $\left[0,\frac kn\right]$ defined by $\phi(y) = H(2y) + 4y + 2(1-2y)\cdot H\(\frac{k/n - y}{1 - 2y}\) - 2H\(\frac kn\)$. Observe that $\phi\(\frac{t_1(n,k)}{n}\) = \psi\(\frac kn\)$.

\noi Let $f$ be a function on $\H$ with $supp\(\widehat{f}\) \subseteq S(n,k)$, such that $R(n,k) = \frac{\E f^4}{\E^2 f^2}$. For an integer $m \ge 1$, consider a function $F_m$ on $nm$ boolean variables defined for $x_1,...,x_m \in \H$ by $F_m\(x_1,...,x_m\) = \prod_{i=1}^m f\(x_i\)$. Observe that for any $\alpha_1,...,\alpha_m \in \H$ holds $\widehat{F_m}\(\alpha_1,...,\alpha_m\) = \prod_{i=1}^m \widehat{f}\(\alpha_i\)$, and hence $supp\(\widehat{F_m}\) \subseteq S(nm,km)$. We also have $\E F_m^p = \(\E f^p\)^m$, for any $p$, and hence $\frac{\E F_m^4}{\E^2 F_m^2} = \(\frac{\E f^4}{\E^2 f^2}\)^m$.
\noi Denoting $N = nm$ and $K = km$, we have that
\[
R(n,k) \queq \frac{\E f^4}{\E^2 f^2} \queq \(\frac{\E F_m^4}{\E^2 F_m^2}\)^{\frac 1m} \quad \le \quad R(N,K)^{\frac 1m}  \quad \le \quad \(C \cdot 2^{\frac{5N}{\log N}} \cdot r(N,K)\)^{\frac 1m}.
\]

\noi Taking $m$ to infinity, we have $R(n,k) \le \liminf_{m \rarrow \infty} \Big(r(N,K)\Big)^{\frac 1m}$. For a sufficiently large $m$ we have $\frac{N}{\log N} \ll K \ll \frac N2 - \frac{N}{\log N}$, and hence, by Lemma~\ref{lem:max-r},
\[
\liminf_{m \rarrow \infty} \Big(r(N,K)\Big)^{\frac 1m} \queq \liminf_{m \rarrow \infty} \Big(\max_t s_t(N,K)\Big)^{\frac 1m} \queq \liminf_{m \rarrow \infty} \Big(\max_{t \in t_1(N,K) \pm \sqrt{N \log N}} s_t(N,K)\Big)^{\frac 1m}
\]
where $t_1(N,K) = \frac{3N - \sqrt{N^2 + 8(N-2K)^2}}{8}$. Recalling that
$s_t(N,K) \queq \frac{1}{{N \choose K}^2} \cdot {N \choose {2t}} \bigg({{2t} \choose t} {{N-2t} \choose {K-t}}\bigg)^2$, and
using the bound ${b \choose a} \le 2^{b H(a/b)}$ (\cite{van-Lint}), we get, for $t \in t_1(N,K) \pm \sqrt{N \log N}$, that
\[
\frac 1n \cdot \log_2 s^{\frac 1m}_t(N,K) \quad \le \quad H\(\frac{2 t}{N}\) + 4 \frac tN + 2 \(1 - 2 \frac tN\) H\(\frac{k/n - t/N}{1 - 2t/N}\) - 2H\(\frac kn\) \queq \phi\(\frac tN\),
\]
where $t/N$ is in $t_1(N,K)/N \pm \sqrt{\frac{\log N}{N}} = t_1(n,k)/n \pm \sqrt{\frac{\log N}{N}}$.
Fixing $n$ and $k$ and taking $m$ to infinity, we get that
\[
\liminf_{m \rarrow \infty} \frac 1n \cdot \log_2 \(\Big(\max_{t \in t_1 \pm \sqrt{N \log N}} s_t(N,K)\Big)^{\frac 1m}\) \quad \le \quad \phi\(\frac{t_1(n,k)}{n}\) \queq \psi\(\frac kn\),
\]
completing the proof of (\ref{thm-claim}).

\noi The remainder of this section is organized as follows. We prove Proposition~\ref{pro:shere-main} in the next subsection. Lemma~\ref{lem:max-r} is proved as one of the steps in that proof. We prove the second inequality of Theorem~\ref{thm:sphere}, namely that
\beqn
\label{thm-second-claim}
2^{n \psi\(\frac kn\)} \quad \le \quad O\(k^{3/2}\) \cdot r(n,k)
\eeqn
in Subsection~\ref{subsec:second claim}. Corollary~\ref{cor:ball} is proved in Subsection~\ref{subsec:ball}.

\subsection{Proof of Proposition~\ref{pro:shere-main}}

\noi We start with observing that by choosing the constant $C$ in the claim of the proposition to be sufficiently large, we may assume that the claim holds for all $n \le n_0$ for any fixed $n_0$ of our choice. Indeed, let $n_0$ be chosen, and set $C = 2^{n_0}$. Then, by the first claim of Proposition~\ref{pro:additive}, for any $n \le n_0$ and $1 \le k \le n/2$ holds
\[
R(n,k) \quad \le \quad {n \choose k} \quad < \quad 2^n \quad \le \quad C \quad \le \quad C \cdot 2^{5 n/\log(n)} \cdot r(n,k)
\]

\noi From now on we fix $n_0$ to be sufficiently large for all asymptotically valid claims below to hold for $n \ge n_0$, and set $C = 2^{n_0}$.

\noi Next, we observe that the claim of the proposition holds when $k$ is very small compared to $n$ or when $k$ is very close to $n/2$. This is done in the next two lemmas.

\lem
\label{lem:small-k}
There exists a sufficiently large constant $n_0$ such that Proposition~\ref{pro:shere-main} holds for all $n \ge n_0$ and $k \le \frac{n}{\log n}$.
\elem
\prf
By (\ref{small-k}) we have that
\[
R(n,k) \quad \le \quad O\(e^{\frac{2k^2}{n}}\) \cdot r(n,k) \quad \le \quad O\(e^{\frac{2n}{\log^2 n}}\) \cdot r(n,k) \quad \le \quad  2^{5 \frac{n}{\log n}} \cdot r(n,k),
\]
for all sufficiently large $n$.

\eprf

\lem
\label{lem:large-k}
There exists a sufficiently large constant $n_0$ such that Proposition~\ref{pro:shere-main} holds for all $n \ge n_0$ and $k \ge \frac n2 - \frac{n}{\log n}$.
\elem
\prf
\noi Assume, w.l.o.g., that $k$ is even. Then, using the inequality ${k \choose {k/2}} \ge \Omega\(\frac{2^k}{\sqrt{k}}\)$, we get
\[
r(n,k) \quad \ge \quad s_{k/2}(n,k) \queq \frac{{n \choose k} \cdot \Big({k \choose {k/2}} {{n-k} \choose {k/2}}\Big)^2} {{n \choose k}^2} \quad > \quad \frac{{k \choose {k/2}}^4}{2^n} \quad \ge \quad \Omega\(\frac{2^{4k - n}}{n^2}\) \quad \ge \quad
\]
\[
\Omega\(2^{n - 4 \frac{n}{\log n} - 2 \log_2 n}\) \quad \ge \quad  2^{n - 5 \frac{n}{\log n}},
\]
where the last inequality holds for a sufficiently large $n$. Therefore,
$R(n,k)~<~2^n~\le~2^{5 \frac{n}{\log n}}~\cdot~r(n,k)$.
\eprf

\noi Hence from now on we may assume that $n$ is sufficiently large and that $\frac{n}{\log n} \le k \le \frac n2 - \frac{n}{\log n}$. The proof of Proposition~\ref{pro:shere-main} will rely on the following two claims.

\pro
\label{pro:sphere-induct}
Let $F(x,y) = \frac{8xy}{4 \sqrt{xy} - \(\sqrt{x} - \sqrt{y}\)^2}$. Then
\begin{enumerate}
\item
The function $F$ is increasing in both $x$ and $y$ in the domain $0 < x/9 < y < 9x$ and is $1$-homogeneous.
\item
For any $1 \le k \le n/2$ the following inductive relation holds: There exist positive numbers $R_0$ and $R_1$ such that $R_0 \le R(n-1,k)$ and $R_1 \le R(n-1,k-1)$ and such that
\[
R(n,k) \quad \le \quad \left\{\begin{array}{ccc} R_0 & \mbox{if} & R_0 \ge 9 R_1 \\ R_1 & \mbox{if} & R_1 \ge 9 R_0 \\ F\(R_0, R_1\) & \mbox{otherwise} \end{array} \right.
\]
\end{enumerate}
\epro

\noi And

\pro
\label{pro:r-F}
There exists a sufficiently large constant $n_0$ such that for all $n \ge n_0$ and for all $\frac{n}{\log n} \le k \le \frac n2 - \frac{n}{\log n}$ holds
\begin{enumerate}

\item
\[
\frac{r(n-1,k-1)}{9} \quad < \quad r(n-1,k) \quad < \quad 9 \cdot r(n-1,k-1)
\]

\item
\[
r(n,k) \quad \in \quad \(1 \pm O\(\frac{\log^{3/2} n}{\sqrt n}\)\) \cdot  F\Big(r(n-1,k-1), r(n-1,k)\Big)
\]
\end{enumerate}

\epro

\noi We first show how to deduce Proposition~\ref{pro:shere-main} from these two claims and then prove the claims.

\noi Assume Proposition~\ref{pro:sphere-induct} and Proposition~\ref{pro:r-F} to hold. Let $n_0$ and $C = 2^{n_0}$ be as defined above. We will argue by induction on $n$ that for all $1 \le k \le n/2$ holds $R(n,k) \le C \cdot 2^{5 \frac{n}{\log n}} \cdot r(n,k)$. Clearly, by the choice of $C$, this holds for $n \le n_0$, which takes care of the base step. We pass to the induction step. Let $1 \le k \le n/2$ be given. We may and will assume that $n \ge n_0$. By Lemmas~\ref{lem:small-k}~and~\ref{lem:large-k} the claim holds for $k \le \frac{n}{\log n}$ and for $k \ge \frac n2 - \frac{n}{\log n}$. So we may assume $\frac{n}{\log n} < k < \frac n2 - \frac{n}{\log n}$.

\noi Let $R_0$ and $R_1$ be the two numbers given by the second claim of Proposition~\ref{pro:sphere-induct}. Consider first the case $R_0 \le \frac{R_1}{9}$. Then, by Proposition~\ref{pro:sphere-induct} and by the induction hypothesis we have
\[
R(n,k) ~ \le ~ R_1 ~ \le ~ R(n-1,k-1) ~ \le ~ C \cdot 2^{5 \frac{n-1}{\log(n-1)}} \cdot r(n-1,k-1) ~ \le ~ C \cdot 2^{5 \frac{n}{\log(n)}} \cdot r(n,k).
\]
Let us explain the last inequality. First, $\frac{n-1}{\log(n-1)} \le \frac{n}{\log(n)}$, since the function $\frac{x}{\log x}$ is increasing for $x \ge e$. Second, simple calculations show that $s_t(n-1,k-1)\leq s_t(n,k) $ for every $0\leq t \leq k$ and $k\leq n/2$, thus $r(n-1,k-1) \le r(n,k)$.

\noi The case $R_1 \le \frac{R_0}{9}$ is treated similarly.

\noi It remains to deal with the case $\frac{R_0}{9} < R_1 < 9 R_0$. In this case, we have $R(n,k) \le F\(R_0,R_1\)$. Let $\rho = \max\left\{\frac{R_0}{r(n-1,k)}, \frac{R_1}{r(n-1,k-1)}\right\}$. Note that by the induction hypothesis $\rho \le C \cdot 2^{5 \frac{n-1}{\log(n-1)}}$. By Proposition~\ref{pro:r-F} the point $\Big(r(n-1,k), r(n-1,k-1)\Big)$ lies in the domain $\{(x,y) : 0 < x/9 < y < 9x\}$ and hence so is the point $\rho \cdot \Big(r(n-1,k), r(n-1,k-1)\Big)$. By the monotonicity of $F$ in this domain and by its $1$-homogeneity, we have
\[
F\(R_0,R_1\) ~\le~ F\Big(\rho \cdot r(n-1,k), \rho \cdot r(n-1,k-1)\Big) ~=~ \rho \cdot F\Big(r(n-1,k), r(n-1,k-1)\Big) ~\le~
\]
\[
C \cdot 2^{5 \frac{n-1}{\log(n-1)}} \cdot F\Big(r(n-1,k), r(n-1,k-1)\Big)
\]

\noi By Proposition~\ref{pro:r-F}, the last expression is at most $C \cdot 2^{5 \frac{n-1}{\log(n-1)}} \cdot \(1 + c \cdot \frac{\log^{3/2} n}{\sqrt n}\) \cdot r(n,k)$, for some absolute constant $c$. Since for large $x$ we have $\(\frac{x}{\log x}\)' \sim \frac{1}{\log^2 x}$, for a sufficiently large $n$ holds
\[
2^{5 \frac{n-1}{\log(n-1)}} \cdot \(1 + c \cdot \frac{\log^{3/2} n}{\sqrt n}\) \quad \le \quad 2^{5 \frac{n-1}{\log(n-1)} + \frac{c}{\ln 2} \cdot \frac{\log^{3/2} n}{\sqrt n}} \quad \le \quad 2^{5 \frac{n}{\log(n)}},
\]
completing the proof of Proposition~\ref{pro:shere-main}.

\subsubsection{Proof of Proposition~\ref{pro:sphere-induct}}

\noi We start with the first claim of the proposition. The function $F(x,y) = F(x,y) = \frac{8xy}{4 \sqrt{xy} - \(\sqrt{x} - \sqrt{y}\)^2}$ is clearly $1$-homogeneous. It's easy to see that it is defined on the domain $0 < x/9 < y < 9x$. A simple computation shows that $\frac{\partial F}{\partial x}$ is proportional to $3 \sqrt{x} - \sqrt{y}$ and therefore is positive in this domain. Hence $F$ increases in $x$. A similar argument shows that $F$ increases in $y$ as well. This completes the proof of the first claim.

\noi We pass to the second claim of the proposition.

\noi Let $f$ be a function on $\H$ with $supp\(\widehat{f}\) \subseteq S(n,k)$, such that $\frac{\E f^4}{\E^2 f^2} = R(n,k)$. Given a function $h$ on $\H$, we can view it as a pair of functions on the two $(n-1)$-dimensional cubes $\{x \in \H, x_n = 0\}$ and $\{x \in \H, x_n = 1\}$. We write this as $h \leftrightarrow \(h_0, h_1\)$. Let $\widehat{f} \leftrightarrow \Big(\widehat{f}_0,\widehat{f}_1\Big)$, and let $g_0$, $g_1$ be functions on the $(n-1)$-dimensional cube such that $g_i = \sum_{\beta} \widehat{f}_i(\beta) W_{\beta}$. It is easy to see that $f \leftrightarrow \(g_0 + g_1, g_0 - g_1\)$. Note that $supp\(\widehat{g_0}\) \subseteq S(n-1,k)$ and $supp\(\widehat{g_1}\) \subseteq S(n-1,k-1)$. We can now define the parameters $R_0$ and $R_1$. Let $R_0 =\frac{\E g_0^4}{\E^2 g_0^2}$. Then $R_0 \le R(n-1,k)$. Similarly, let $R_1 =\frac{\E g_1^4}{\E^2 g_1^2}$. Then $R_1 \le R(n-1,k-1)$.

\noi A simple calculation and an application of the Cauchy-Schwarz inequality gives
\[
R(n,k) \queq \frac{\E f^4}{\E^2 f^2} \queq \frac{\E g_0^4 + 6 \E g_0^2 g_1^2 + \E g_1^4}{\E^2 g_0^2 + 2 \E g_0^2 \E g_1^2 + \E^2 g_1^2} \quad \le \quad \frac{\E g_0^4 + 6 \sqrt{\E g_0^4 \E g_1^4} + \E g_1^4}{\E^2 g_0^2 + 2 \E g_0^2 \E g_1^2 + \E^2 g_1^2}
\]

\noi Let $m\(g_0,g_1\)$ be the supremum of the RHS over a $1$-parameter family of expressions, where we replace $g_1$ with $\theta \cdot g_1$, for a real parameter $\theta$. Clearly $R(n,k) \le m\(g_0,g_1\)$.
We will show that
\[
m\(g_0,g_1\) \queq \left\{\begin{array}{ccc} R_0 & \mbox{if} & R_0 > 9 R_1 \\ R_1 & \mbox{if} & R_1 > 9 R_0 \\ F\(R_0, R_1\) & \mbox{otherwise,} \end{array} \right.
\]
and this will complete the proof of the proposition.

\noi Consider the following function of a nonnegative parameter $x = \theta^2$:
\[
G(x) \queq \frac{\E g_1^4 \cdot x^2  + 6 \sqrt{\E g_0^4 \E g_1^4} \cdot x  + \E g_0^4}{\E^2 g_1^2 \cdot x^2  + 2 \E g_0^2 \E g_1^2 \cdot x + \E^2 g_0^2} 
\]

\noi By definition $m\(g_0,g_1\) = \mbox{sup}_{x \ge 0} G(x)$. It is easy to see that the derivative $G'$ equals, up to a positive factor, to
\[
Q(x) \queq \sqrt{\E g_1^4} \cdot \(\sqrt{R_1} - 3 \sqrt{R_0}\) \cdot x ~+~ \sqrt{\E g_0^4} \cdot \(3\sqrt{R_1} -  \sqrt{R_0}\).
\]

\noi If $R_0 \ge 9 R_1$, then $Q \le 0$, which means $G$ is decreasing on $(0,\infty)$, and $m\(g_0,g_1\) = G(0) = R_0$. Similarly, if $R_1 \ge 9 R_0$, then $G$ increases on $(0,\infty)$, and $m\(g_0,g_1\) = G(\infty) = R_1$.

\noi The interesting case is when $1/9 R_1 < R_0 < 9 R_1$, and then the unique maximum of $G$ is attained at the root of $Q$, that is at
\[
x \queq \frac{\sqrt{\E g_0^4}}{\sqrt{\E g_1^4}} \cdot \frac{3\sqrt{R_1} -  \sqrt{R_0}}{3\sqrt{R_0} -  \sqrt{R_1}}.
\]

\noi Substituting this value of $x$ and simplifying we get that
\[
m\(g_0,g_1\) \queq G(x) \queq \frac{8R_0R_1}{4 \sqrt{R_0 R_1} - \(\sqrt{R_0} - \sqrt{R_1}\)^2} \queq F\(R_0,R_1\),
\]
completing the proof of the proposition.

\subsubsection{Proof of Proposition~\ref{pro:r-F}}
\noi From now on we assume (in this subsection) that the assumptions of the proposition hold, that is that
$n$ is sufficiently large and that $\frac{n}{\log n} \le k \le \frac n2 - \frac{n}{\log n}$.

\noi Recall that $r(n,k) = \frac{1}{{n \choose k}^2} \cdot \sum_{t = 0}^k {n \choose {2t}} \({{2t} \choose t} \cdot {{n-2t} \choose {k-t}}\)^2 = \sum_{t=0}^k s_t(n,k)$, where
$s_t(n,k) = \frac{1}{{n \choose k}^2} \cdot {n \choose {2t}} \({{2t} \choose t} \cdot {{n-2t} \choose {k-t}}\)^2$. We start with the following claim.
\lem
\label{lem:st-main}
Let $t_1(n,k) =  \frac{3n - \sqrt{n^2 + 8(n-2k)^2}}{8}$. Then
\begin{enumerate}

\item $\frac k3 \le t_1(n,k) \le \frac{11k}{12}$.

\item
Let $3 \le \Delta < t_1(n,k)$. Then for $1 \le t \le t_1(n,k) - \Delta$ holds
\[
\frac{s_{t+1}}{s_t} \quad \ge \quad 1 + \frac{\Delta}{t}
\]

\item
Let $\log n \le \Delta < k - t_1(n,k)$. Then for $t_1(n,k) + \Delta \le t \le k-1$ holds
\[
\frac{s_{t+1}}{s_t} \quad \le \quad 1 - \frac{\Delta}{t}
\]
\end{enumerate}

\elem

\noi We record two immediate corollaries of this lemma. Choosing $\Delta = \sqrt{n \log n}$, we obtain
\[
\max_{0 \le t \le k} s_t(n,k) \queq \max_{t \in t_1(n,k) \pm \sqrt{n \log n}} ~s_t(n,k),
\]
which is the claim of Lemma~\ref{lem:max-r}. Another immediate corollary is that

\cor
\label{cor:central interval}
There is an interval of length $L = L(n) = O\(\sqrt{n \log n}\)$ such that
\[
r(n,k) \quad \le \quad \(1 + \frac 1n\) \cdot \sum_{t_1 - L \le t \le t_1 + L} ~s_t
\]
\ecor

\prf (of Lemma~\ref{lem:st-main})

\noi We start with the first claim of the lemma. We have that
\[
t_1(n,k) ~=~ \frac{3n - \sqrt{n^2 + 8(n-2k)^2}}{8} ~=~ \frac{32k(n-k)}{8 \cdot (3n + \sqrt{n^2 + 8(n-2k)^2})} ~\ge~ \frac{4k(n-k)}{6n} ~\ge~ \frac k3,
\]
where the last inequality holds since $k \le n/2$.

\noi On the other hand, $k - t_1(n,k) = \frac{\sqrt{n^2 + 8(n-2k)^2} - (3n - 8k)}{8}$. If $k \ge 3n/8$, this is at least $\frac{\sqrt{n^2 + 8(n-2k)^2}}{8} \ge \frac n8 \ge \frac k4$. Otherwise, if $k < 3n/8$, this equals
\[
\frac{n^2 + 8(n-2k)^2 - (3n-8k)^2}{8\cdot\(\sqrt{n^2 + 8(n-2k)^2}  + (3n - 8k)\)} ~=~ \frac{2k(n-2k)}{\sqrt{n^2 + 8(n-2k)^2}  + (3n - 8k)} ~>~
\frac{2k \cdot n/4}{6n} ~\ge~ \quad \frac{k}{12}.
\]

\noi We proceed to the second and the third claims. Consider the ratio $s_{t+1} / s_t$. After some simplifying, this ratio is
\[
\frac{{n \choose {2t+2}} \({{2t+2} \choose {t+1}} {{n-2t-2} \choose {k-t-1}}\)^2}{{n \choose {2t}} \({{2t} \choose t} {{n-2t} \choose {k-t}}\)^2} \queq
\frac{2 (2t+1)}{(t+1)^3} \cdot \frac{(k-t)^2 (n-k-t)^2}{(n-2t) (n-2t-1)} \queq
\]
\[
\(\frac{2(k-t)(n-k-t)}{t(n-2t)}\)^2 \cdot \Big(1 + \e(n,t)\Big), \quad \mbox{where} \quad 1 + \e(n,t) ~=~ \frac{(2t+1)t^2}{2(t+1)^3} \cdot \frac{n-2t}{n-2t-1}
\]
\noi It is easy to see that, in our assumptions for $k$ and $n$, we have $-\frac{3}{t} \le \e(n,t) \le \frac{1}{n-2t}$. We introduce some notation. Let $r(t)= \frac{2(k-t)(n-k-t)}{t(n-2t)}$, and let $q(t) = 4t^2 - 3nt +2k(n-k)$. Then $\frac{s_{t+1}}{s_t} = r^2(t) \cdot \big(1 + \e(n,t)\big)$, and $r(t) = 1 + \frac{q(t)}{t(n-2t)}$.

\noi  The roots of the quadratic $q(t)$ are $t_{1,2}(n,k) = \frac{3n \pm \sqrt{n^2 + 8(n-2k)^2}}{8}$. (From now on till the end of this subsection we write $t_1,t_2$ for $t_1(n,k)$ and $t_2(n,k)$.) We know that $t_1 < k$ and it is easy to see that $t_2 > n/2 > k$. Hence, for $t \le t_1 - \Delta$ we have:
\[
r(t) \queq 1 + \frac{q(t)}{t(n-2t)} \queq  1 + \frac{4\(t-t_1\)\(t-t_2\)}{t(n-2t)} \quad \ge \quad 1 + \frac{4\Delta \cdot (t_2 - t_1)}{t(n-2t)} \queq
\]
\[
\frac{\Delta \cdot \sqrt{n^2 + 8(n-2k)^2}}{t(n-2t)} \quad \ge \quad
1 + \frac{\Delta n}{t(n-2t)} \quad \ge \quad 1 + \frac{\Delta}{t}
\]
Therefore, $\frac{s_{t+1}}{s_t} = r^2(t) \cdot \big(1 + \e(n,t)\big) \ge 1 + \frac{2\Delta}{t} - \frac{3}{t} \ge 1 + \frac{\Delta}{t}$. This completes the proof of the second claim of the lemma.

\noi For $t_1 + \Delta \le t \le k-1$, we have
\[
r(t) \queq 1 + \frac{4\(t-t_1\)\(t-t_2\)}{t(n-2t)} \quad \le \quad  1 - \frac{4\Delta \cdot \(t_2 - t\)}{t(n-2t)} \quad \le \quad 1 - \frac{4\Delta \cdot \(\frac n2 - t\)}{t(n-2t)} \queq 1 - \frac{2\Delta}{t}
\]

\noi Therefore $\frac{s_{t+1}}{s_t} \le 1 - \frac{2\Delta}{t} + \frac{1}{n-2t} \le 1 - \frac{\Delta}{t}$, recalling that $\Delta \ge \log(n)$, and $n - 2t \ge n - 2k \ge 2\frac{n}{\log n}$.

\eprf

\cor
\label{cor:ratios}
\begin{enumerate}
\item
\[
\frac{r(n,k-1)}{r(n,k)} \in \(1 \pm O\(\frac{\log^{3/2} n}{\sqrt n}\)\) \cdot \(\frac{n-k}{k} \cdot \frac{k-t_1}{n-k-t_1}\)^2
\]

\item
\[
\frac{r(n-1,k-1)}{r(n,k)} \in \(1 \pm O\(\frac{\log^{3/2} n}{\sqrt n}\)\) \cdot \frac{n}{n-2t_1} \cdot \(\frac{k-t_1}{k}\)^2
\]

\end{enumerate}
\ecor

\prf
We prove only the first claim of the corollary. The proof of the remaining claim is similar.

\noi Note that $| t_1(n,k) - t_1(n,k-1) | \le 1$. Let $L = O\(\sqrt{n \log n}\)$ be the length of the interval around $t_1 = t_1(n,k)$ such that both $r(n,k)$ and $r(n,k-1)$ are attained, up to an $\(1 - 1/n\)$-factor by summing the corresponding summands in this interval (by Corollary~\ref{cor:central interval}). It suffices to show that for any $t$ in the interval $t_1 \pm L$ holds $s_t(n,k-1) / s_t(n,k) \in 1 \pm O\(\frac{\log^{3/2} n}{\sqrt n}\)$. Indeed we have
\[
\frac{s_t(n,k-1)}{s_t(n,k)} \queq \(\frac{n-k}{k} \cdot \frac{k-t}{n-k-t}\)^2 \queq
\]
\[
\(\frac{k-t}{k-t_1} \cdot \frac{n-k-t_1}{n-k-t}\)^2 \cdot \(\frac{n-k}{k} \cdot \frac{k-t_1}{n-k-t_1}\)^2 \in
\]
\[
\(1 \pm O\(\frac{L}{k}\)\) \cdot \(1 \pm O\(\frac{L}{n}\)\) \cdot \(\frac{n-k}{k} \cdot \frac{k-t_1}{n-k-t_1}\)^2 \subseteq
\]
\[
\(1 \pm O\(\frac{\log^{3/2} n}{\sqrt n}\)\) \cdot \(\frac{n-k}{k} \cdot \frac{k-t_1}{n-k-t_1}\)^2
\]

\eprf

\noi We are now ready to prove Proposition~\ref{pro:r-F}. The first claim of the proposition is that $1/9 <  \frac{r(n,k-1)}{r(n,k)} < 9$. In fact, it is easy to see that for all $0 \le t \le k-1$ holds $\frac{s_t(n,k-1)}{s_t(n,k)} = \(\frac{n-k}{k} \cdot \frac{k-t}{n-k-t}\)^2 \le 1$, so the upper bound trivially holds, even with $9$ replaced by $1$. We pass to the lower bound. By the first claim of Corollary~\ref{cor:ratios} it suffices to show that for some absolute constant $c > 0$ holds
$1/3 + c/\log(n) \le \frac{n-k}{k} \cdot \frac{k-t_1}{n-k-t_1}$. We write $\delta$ for $c/\log(n)$.

\noi After rearranging, we need to show that $t_1 \le \frac{\(2/3 - \delta\)k(n-k)}{(n-k) - \(1/3+\delta\)k}$. This would follow from a stronger inequality $t_1 \le \(1-\frac{3\delta}{2}\) \cdot \frac{2k(n-k)}{3n-4k}$. Recall that $t_1$ is a root of the quadratic $4t^2 - 3nt + 2k(n-k)$. Substituting $3nt_1 - 4t^2_1$ for $2k(n-k)$ it is easy see that this inequality would follow from $\frac{3n - 4t_1}{3n - 4k} \ge 1 + 3 \delta$. Recall that $t_1 \le \frac{11k}{12}$, and that $k \ge \frac{n}{\log(n)}$. Hence $\frac{3n - 4t}{3n - 4k} \ge 1 + \frac{k}{3(3n-4k)} \ge 1 + \frac{1}{9 \log(n)}$, completing the proof (for $c$ small enough).
\eprf

\noi We pass to the second claim of the proposition. Let $x = r(n-1,k-1)$. Let $y = \frac{k^2\(n-k-t_1\)^2}{(n-k)^2 \(k-t_1\)^2} \cdot x$. Let $z = \frac{\(n-2t_1\)k^2}{n \(k-t_1\)^2} \cdot x$. By Corollary~\ref{cor:ratios} we have that $y \in \(1 \pm O\(\frac{\log^{3/2} n}{\sqrt n}\)\) \cdot r(n-1,k)$ and $z \in \(1 \pm O\(\frac{\log^{3/2} n}{\sqrt n}\)\) \cdot r(n,k)$.

\noi Next we claim that $z = F(x,y)$. Since $F$ is $1$-homogeneous, it suffices to verify the identity
\[
\frac{\(n-2t_1\)k^2}{n \(k-t_1\)^2} \queq F\(1,\frac{k^2\(n-k-t_1\)^2}{(n-k)^2 \(k-t_1\)^2}\).
\]

\noi Simplifying, it is the same as showing:
\[
\frac{n-2t_1}{n} \queq \frac{8\(n-k-t_1\)^2\(k-t_1\)^2}{6k(n-k)\(n-k-t_1\)\(k-t_1\) - k^2\(n-k-t_1\)^2 - (n-k)^2\(k-t_1\)^2}.
\]
This can be verified by applying several times the identity $4t_1^2 - 3nt_1 + 2k(n-k) = 0$. We omit the details.

\noi Now we can conclude the proof. Let $\rho = \max\left\{\frac{r(n-1,k)}{y}, 1\right\}$. Then $\rho \le 1 + O\(\frac{\log^{3/2} n}{\sqrt{n}}\)$. By the proof of the first claim of the proposition, the point $(x,y)$ lies in the domain $0 < x/9 < y < 9x$ and hence also the point $(\rho \cdot x, \rho \cdot y)$. Both coordinates of this point are larger or equal to those of $\Big(r(n-1,k-1), r(n-1,k)\Big)$, which, by the first claim of the proposition, also lies in this domain. By the $1$-homogeneity and monotonicity of $F$ in this domain we have
\[
F\Big(r(n-1,k-1), r(n-1,k)\Big) \le  F(\rho x, \rho y) = \rho F(x,y) = \rho z \in \(1 \pm O\(\frac{\log^{3/2} n}{\sqrt n}\)\) \cdot r(n,k).
\]
\eprf

\subsection{Proof of (\ref{thm-second-claim})}
\label{subsec:second claim}

\noi We start with some simple observations. First, as above, by making the constant hidden in the asymptotic notation to be large enough, we may assume that the claim holds for $n \le n_0$, for any fixed $n_0$ of our choice. Next, it suffices to show the claim for $k \ge k_0$, for any fixed $k_0$ that we choose. This is because for $k < k_0$ the set $S(n,k)$ may be viewed as a subset of $S\big(n + (k_0 - k), k_0\big)$ (see a similar argument in the proof of Proposition~\ref{pro:shere-main}). From now on we assume $n \ge n_0$ and $k \ge k_0$, for sufficiently large $n_0$ and $k_0$.

\noi We will work with the function $\phi$ introduced in the proof of Theorem~\ref{thm:sphere}. Recall that $\phi$ is a function on $\left[0,\frac kn\right]$ defined by $\phi(y) = H(2y) + 4y + 2(1-2y)\cdot H\(\frac{k/n - y}{1 - 2y}\) - 2H\(\frac kn\)$, and that $\psi\(\frac kn\) = \phi\(\frac{t_1(n,k)}{n}\)$.

\noi Let $t_1$ stand for $t_1(n,k)$, and let $t^{\ast} = \lceil t_1 \rceil$. We proceed as follows:  First, we observe that $\phi$ is defined on $\frac{t^{\ast}}{n}$ and that $2^{n \phi\(\frac{t^{\ast}}{n}\)} \le  O\(k^{3/2}\) \cdot r(n,k)$. Then we show that $\phi\(\frac{t_1}{n}\)$ and $\phi\(\frac{t^{\ast}}{n}\)$ differ by at most $O\(\frac 1n\)$, which implies $2^{n \psi\(\frac{k}{n}\)} = 2^{n \phi\(\frac{t_1}{n}\)} \le O\(2^{n \phi\(\frac{t^{\ast}}{n}\)}\)$, and completes the proof.

\noi By Lemma~\ref{lem:st-main}, $t_1 \le \frac{11 k}{12}$, and hence $t^{\ast} \le t_1 + 1 \le k$. Therefore $\phi$ is defined on $\frac{t^{\ast}}{n}$. Next, recall that, by Stirling's formula, for all $0 < a < b$ holds ${b \choose a} = \Theta\(\sqrt{\frac{b}{a(b-a)}}\) \cdot~2^{b H\(a/b\)}$, and in particular, for $0 < a \le b/2$ holds ${b \choose a} = \Theta\(\sqrt{\frac{1}{a}}\) \cdot 2^{b H\(a/b\)}$. Substituting this estimate for the binomial coefficients in the formula for $s_{t^{\ast}}(n,k)$ gives
\[
r(n,k) ~\ge~ s_{t^{\ast}}(n,k) ~=~ \Theta\(\frac{k}{t^{\ast}\(k-t^{\ast}\)} \cdot \sqrt{\frac{n}{2t^{\ast}\(n-2t^{\ast}\)}}\) \cdot 2^{n \phi\(\frac{t^{\ast}}{n}\)} ~\ge~ \Omega\(k^{-\frac32}\) \cdot 2^{n \phi\(\frac{t^{\ast}}{n}\)}.
\]

\noi Next, we argue that $\Big |\phi\(\frac{t_1}{n}\) - \phi\(\frac{t^{\ast}}{n}\) \Big | \le O\(\frac 1n\)$. Since $t_1 \le t^{\ast} < t_1 + 1$, it suffices to show that the absolute value of the derivative of $\phi$ is bounded by a constant on $\(\frac{t_1}{n}, \frac{t^{\ast}}{n}\)$. Let $a := \frac kn$. Then
\beqn
\label{der-phi}
\frac12 \phi'(y) \queq \log\(\frac{1-2y}{2y}\) + 2 - 2H\(\frac{a-y}{1-2y}\) - \frac{1-2a}{1-2y} \cdot \log\(\frac{1-a-y}{a-y}\)
\eeqn
Let $\frac{t_1}{n} < y < \frac{t^{\ast}}{n}$. Then, by Lemma~\ref{lem:st-main} and by our assumptions on $k$ and $n$, we have that $0 < a \le \frac12$ and $c_1 a \le y \le \(1-c_2\) a$, for some absolute constants $0 < c_1, c_2 < 1$. It is easy to see that for $a$ bounded away from zero all the terms on the RHS of (\ref{der-phi}) are bounded. Hence it only remains to consider the case $a \rarrow 0$. To deal with this case, we can rewrite (\ref{der-phi}) as follows (omitting the second and the third term on the RHS, since their contribution is bounded by $2$):
\[
\frac12 \phi'(y) \quad \approx \quad  \log\(\frac{1-2y}{2y}\) - \frac{1-2a}{1-2y} \cdot \log\(\frac{1-a-y}{a-y}\) \queq
\]
\[
\(\log\(\frac{1-2y}{2y}\) - \log\(\frac{1-a-y}{a-y}\)\) + 2\frac{a-y}{1-2y} \cdot \log\(\frac{1-a-y}{a-y}\) \queq
\]
\[
\log\(\frac{1-2y}{1-a-y}\) + \log\(\frac{a-y}{2y}\) + 2\frac{a-y}{1-2y} \cdot \log\(\frac{1-a-y}{a-y}\).
\]
It is easy to see that all the summands in the last expression are bounded by a constant, completing the proof of (\ref{thm-second-claim}).

\eprf

\subsection{Proof of Corollary~\ref{cor:ball}}
\label{subsec:ball}

\noi We start with the first claim. Let $A \subseteq \H$ be a Hamming ball of radius $k$. First, we observe that Proposition~\ref{pro:shere-main} implies the following bound on $\mu(A)$.
\beqn
\label{weak-ball}
\mu(A) \quad \le \quad C \cdot  (k+1)^3  2^{5 n/\log(n)} \cdot r(n,k).
\eeqn

\noi To see this, let $f$ be a function on $\H$ with $supp\(\widehat{f}\) \subseteq A$. Write $f = \sum_{i=0}^k f_i$, with $supp\(\widehat{f}\) \subseteq S(n,i)$, for $i = 0,...,k$. By Proposition~\ref{pro:shere-main}, we have $\E f^4_i \le C \cdot 2^{5 n/\log(n)} \cdot r(n,i) \cdot \E^2 f^2_i$. In the proof of Proposition~\ref{pro:r-F}, we have observed that $s_t(n,k-1) \le s_t(n,k)$, for all $0 \le t \le k-1$, which implies $r(n,k-1) \le r(n,k)$, and hence $r(n,i) \le r(n,k)$, for all $0 \le i \le k$. Consequently, we have $\E f^4_i \le C \cdot  2^{5 n/\log(n)} \cdot r(n,k) \cdot \E^2 f^2_i$. Observing that the functions $\{f_i\}$ are orthogonal, and using Jensen's inequality, we have:
\[
\E f^4 ~\le~ (k+1)^3 \cdot \sum_{i=0}^k \E f^4_i  ~\le~  C \cdot (k+1)^3  2^{5 n/\log(n)} \cdot r(n,k) \cdot \sum_{i=0}^k {\E}^2 f^2_i  ~ \le ~
\]
\[
 C \cdot (k+1)^3  2^{5 n/\log(n)} \cdot r(n,k) \cdot \(\sum_{i=0}^k \E f^2_i\)^2 ~=~ C \cdot (k+1)^3  2^{5 n/\log(n)} \cdot r(n,k) \cdot {\E}^2 f^2,
\]
completing the proof of (\ref{weak-ball}).

\noi The inequality $\mu(A) \le  2^{n \psi\(\frac kn\)}$ can now be derived from (\ref{weak-ball}) by a 'tensorization argument', as in derivation of the first claim of Theorem~\ref{thm:sphere} from  Proposition~\ref{pro:shere-main}. We omit the details.

\noi We pass to the second claim. It suffices to show that $\psi(x) \le \min\left\{2\log_2(3) \cdot x, 1\right\}$ for all $0 \le x \le 1/2$, and moreover $\psi(x) = 2\log_2(3) \cdot x$ only at $x = 0$, and $\psi(x) = 1$ only at $x = 1/2$. The key observation is that $\psi$ is strongly concave.

\lem
\label{lem:psi}
For all $0 < x \le 0.5$ holds $\psi''(x) < 0$.
\elem

\prf

\noi We have that
\[
\psi'(x) \queq 2r' \cdot \log_2\(\frac{1-2r}{2r}\) + 4r' \cdot \(1 - H\(\frac{x-r}{1-2r}\)\) +
\]
\[
2\cdot\(1 - \frac{(1-2x) r'}{1-2r}\) \cdot \log_2\(\frac{1-x-r}{x-r}\) - 2\log_2\(\frac{1-x}{x}\),
\]
and, after some rearrangement, that
\[
\frac12 \cdot \psi''(x) \queq r'' \cdot \Bigg(\log_2\(\frac{1-2r}{2r}\) + 2 \cdot \(1 - H\(\frac{x-r}{1-2r}\)\) - \frac{1-2x}{1-2r} \cdot \log_2\(\frac{1-x-r}{x-r}\)\Bigg) ~-~
\]
\[
\frac{\(r'\)^2}{\ln 2 \cdot r(1-2r)} ~-~ \frac{\big((1-2r) - (1-2x)r'\big)^2}{\ln 2 \cdot (1-2r)(x-r)(1-x-r)} ~+~ \frac{1}{\ln 2 \cdot x(1-x)}
\]

\noi We claim that the term which multiplies $r''$ is zero. To see that, we make some observations about the function $r$, which will be useful later on as well. First, it is easy to see that it increases from $0$ to $0.25$ on $[0,0.5]$. Next, we have $r' = \frac{2-4x}{3-8r}$, and finally the identity $\frac12 \(3r - 4r^2\) = x(1-x)$, which follows e.g., from the fact that $t_1(n,k)$ is a root of the quadratic $4t^2 - 3nt +2k(n-k) = 0$.

\noi Next, after some simplifying, we have
\[
\log_2\(\frac{1-2r}{2r}\) + 2 \cdot \(1 - H\(\frac{x-r}{1-2r}\)\) - \frac{1-2x}{1-2r} \cdot \log_2\(\frac{1-x-r}{x-r}\) = \log_2\(\frac{2(x-r)(1-x-r)}{r(1-2r)}\).
\]

\noi Using the identity $\frac12 \cdot \(3r - 4r^2\) = x(1-x)$, it is easy to see that $2(x-r)(1-x-r) = r(1-2r)$ and hence the RHS vanishes. This simplifies the expression for $\psi''$ to:
\[
\psi''(x) \queq -\frac{2}{\ln 2} \cdot \(\frac{\(r'\)^2}{r(1-2r)} + \frac{\big((1-2r) - (1-2x)r'\big)^2}{(1-2r)(x-r)(1-x-r)} - \frac{1}{x(1-x)}\)
\]

\noi Since $r' = \frac{2-4x}{3-8r}$, we have $\(r'\)^2 = \frac{4(1-2x)^2}{(3-8r)^2} = \frac{4(1-2r)(1-4r)}{(3-8r)^2}$. Similarly, $(1-2x)r' = \frac{2(1-2x)^2}{3-8r} = \frac{2(1-2r)(1-4r)}{3-8r}$. Making these substitutions, replacing $(x-r)(1-x-r)$ with $\frac12 r(1-2r)$ and $x(1-x)$ with $\frac12 \(3r - 4r^2\)$, and simplifying, we get
\[
\frac{\(r'\)^2}{r(1-2r)} + \frac{\big((1-2r) - (1-2x)r'\big)^2}{(1-2r)(x-r)(1-x-r)} - \frac{1}{x(1-x)} \queq \frac{8}{(3-8r)(3-4r)} \quad > \quad 0,
\]
completing the proof of the lemma.

\eprf

\noi We can now complete the proof of the second claim of the corollary. It is easy to see that $\psi'\(1/2\) = 0$. Since $\psi''$ is negative, this means that $\psi'$ is positive on $\(0,1/2\)$ and hence the unique maximum of $\psi$ is at $1/2$, where it equals $1$.

\noi On the other hand, using the fact that $r'(0) = 2/3$, it is easy to see that $\lim_{x \rarrow 0} \psi'(x) = 2 \log_2(3)$. Since $\psi''$ is negative, this means that $\psi' < 2 \log_2(3)$ on $\(0,1/2\)$ and hence that $\psi(x) < 2 \log_2(3) \cdot x$ for all $0 < x \le 1/2$.

\eprf

\subsubsection*{Acknowledgement}

\noi We are grateful to Yury Polyanskiy for many valuable remarks, in particular for pointing out that  our results are relevant to the questions investigated in \cite{P2}. We also thank Yuzhou Gu for valuable remarks.


\begin{thebibliography}{99}

\bibitem{Bonami}
A. Bonami, {\sl Etude des coefficients Fourier des fonctions de Lp(G)}, Annales
de l’Institut Fourier, 20(2) (1970), 335–402.


\bibitem{Don}
D. L. Donoho and P. B. Stark, {\sl Uncertainty principles and signal recovery}, SIAM J. Applied
Math., 49(1989), 906-931.

\bibitem{BGS}
W. T. Gowers, {\sl A new proof of Szemeredi's theorem for arithmetic progressions of length four},
GAFA 8 (1998), 529-551.

\bibitem{GT}
B. Green and T. Tao, {\sl Freiman's Theorem in Finite Fields via Extremal Set Theory},  Combinatorics, Probability \& Computing, Vol. 18(3) (2009), 335-355.


\bibitem{Hastad-personal}
J. Hastad, personal communication.

\bibitem{KM}
J. Kahn and R. Meshulam, {\sl Uncertainty inequalities on Hamming cubes}, unpublished (1996).

\bibitem{van-Lint}
J. H. van Lint, {\bf Introduction to coding theory}, Springer-Verlag, Berlin, 1999.

\bibitem{NO}
P. Nayar and K. Oleszkiewicz, {\sl Khinchine type inequalities with optimal
constants via ultra log-concavity}, Positivity, 2012.

\bibitem{O'Donnel}
R. O'Donnel, {\bf Analysis of Boolean functions}, Cambridge University Press, 2014.

\bibitem{P2}
Y. Polyanskiy, {\sl Hypercontractivity of spherical averages in Hamming space}, arXiv:1309.3014, 2013.

\bibitem{PS}
Y. Polyanskiy and A. Samorodnitsky, {\sl Improved log-Sobolev inequalities, hypercontractivity and uncertainty principle on the hypercube}, arXiv:1606.07491, 2016.

\bibitem{Sanders}
T. Sanders, {\sl On the Bogolyubov-Ruzsa lemma}, Analysis \& PDE, 5(3), 2012, 627-655.

\bibitem{Shkredov}
I.D. Shkredov, {\sl An introduction to higher energies and sumsets}, arXiv:1512.00627, 2015.

\bibitem{TV}
T. Tao and V. Vu, {\bf Additive Combinatorics}, Cambridge University Press 2006.

\end{thebibliography}
\end{document}